\documentclass[a4paper,12pt]{amsart}
\usepackage{amsfonts}
\usepackage{amssymb}
\usepackage{ifthen}
\usepackage{graphicx}
\nonstopmode \numberwithin{equation}{section}
\newtheorem{thm}{Theorem}
\newtheorem{cor}{Corollary}
\newtheorem{lem}{Lemma}
\newtheorem{prop}{Proposition}

\newtheorem{conj}[equation]{Conjecture}

\theoremstyle{definition}
\newtheorem{defn}{Definition}

\newtheorem{prob}[equation]{Problem}
\newtheorem{ques}[equation]{Question}
\newtheorem{rem}{Remark}[section]

\newcounter {own}
\def\theown {\thesection       .\arabic{own}}

\newenvironment{pf}[1][]{%
 \vskip 2mm
 \noindent
 \ifthenelse{\equal{#1}{}}%
  {{\slshape Proof. }}%
  {{\slshape #1.} }%
 }%
{\qed\bigskip}

\newcounter{alphabet}
\newcounter{tmp}

\makeatletter
\newcommand{\Ref}[1]{\@ifundefined{r@#1}{}{\setcounter{tmp}{\ref{#1}}\Alph{tmp}}}
\makeatother

\newenvironment{Lem}[1][]{\refstepcounter{alphabet}%
\bigskip%
\noindent%
{\bf Lemma \Alph{alphabet}}%
{\bf .} \itshape}{\vskip 6pt}

\newcommand{\ID}{{\mathbb D}}

\def\be{\begin{equation}}
\def\ee{\end{equation}}

\newcommand{\bee}{\begin{enumerate}}
\newcommand{\eee}{\end{enumerate}}

\newcommand{\blem}{\begin{lem}}
\newcommand{\elem}{\end{lem}}
\newcommand{\bthm}{\begin{thm}}
\newcommand{\ethm}{\end{thm}}
\newcommand{\bcor}{\begin{cor}}
\newcommand{\ecor}{\end{cor}}
\newcommand{\beg}{\begin{examp}}
\newcommand{\eeg}{\end{examp}}
\newcommand{\begs}{\begin{examples}}
\newcommand{\eegs}{\end{examples}}
\newcommand{\bdefe}{\begin{defn}}
\newcommand{\edefe}{\end{defn}}
\newcommand{\bprob}{\begin{prob}}
\newcommand{\eprob}{\end{prob}}
\newcommand{\bques}{\begin{ques}}
\newcommand{\eques}{\end{ques}}
\newcommand{\bei}{\begin{itemize}}
\newcommand{\eei}{\end{itemize}}

\newcommand{\bcon}{\begin{conj}}
\newcommand{\econ}{\end{conj}}
\newcommand{\bcons}{\begin{conjs}}
\newcommand{\econs}{\end{conjs}}
\newcommand{\bprop}{\begin{prop}}
\newcommand{\eprop}{\end{prop}}
\newcommand{\br}{\begin{rem}}
\newcommand{\er}{\end{rem}}
\newcommand{\brs}{\begin{rems}}
\newcommand{\ers}{\end{rems}}
\newcommand{\bo}{\begin{obser}}
\newcommand{\eo}{\end{obser}}
\newcommand{\bos}{\begin{obsers}}
\newcommand{\eos}{\end{obsers}}
\newcommand{\bpf}{\begin{pf}}
\newcommand{\epf}{\end{pf}}
\newcommand{\ba}{\begin{array}}
\newcommand{\ea}{\end{array}}
\newcommand{\beq}{\begin{eqnarray}}
\newcommand{\beqq}{\begin{eqnarray*}}
\newcommand{\eeq}{\end{eqnarray}}
\newcommand{\eeqq}{\end{eqnarray*}}

\newcommand{\ds}{\displaystyle}

\begin{document}
\bibliographystyle{amsplain}
\title [On some properties of solutions of the $p$-harmonic
equation]{On properties of solutions of the $p$-harmonic equation}
\author{SH. Chen}
\address{SH. Chen, Department of Mathematics,
Hunan Normal University, Changsha, Hunan 410081, People's Republic
of China.} \email{shlchen1982@yahoo.com.cn}
\author{S. Ponnusamy}
\address{S. Ponnusamy, Department of Mathematics,
Indian Institute of Technology Madras, Chennai-600 036, India.}
\email{samy@iitm.ac.in}
\author{X. Wang${}^{~\mathbf{*}}$}
\address{X. Wang, Department of Mathematics,
Hunan Normal University, Changsha, Hunan 410081, People's Republic
of China.} \email{xtwang@hunnu.edu.cn}

\subjclass[2000]{Primary: 30C65, 30C45; Secondary: 30C20}
\keywords{$p$-harmonic mapping, starlikeness, convexity, region of variability,
 Landau's theorem\\
${}^{\mathbf{*}}$ Corresponding author}

\begin{abstract}
A $2p$-times continuously differentiable complex-valued function
$f=u+iv$ in a simply connected domain $\Omega\subseteq\mathbb{C}$ is \textit{p-harmonic}
if $f$ satisfies the $p$-harmonic equation
$\Delta ^pf=0.$
In this paper, we investigate the properties of $p$-harmonic
mappings in the unit disk $|z|<1$. First, we discuss the convexity,
the starlikeness and the region of variability of some classes of
$p$-harmonic mappings. Then we prove the existence of Landau
constant for the class of functions of the form $Df=zf_{z}-\overline
zf_{\overline z}$, where $f$ is $p$-harmonic in $|z|<1$.
Also, we discuss the region  of variability for certain $p$-harmonic mappings.
At the end, as a consequence of the earlier results of the authors,
we present explicit upper estimates for Bloch norm for bi- and tri-harmonic mappings.
\end{abstract}

\thanks{The research was partly supported by NSFs of China (No. 11071063)
}

\maketitle \pagestyle{myheadings} \markboth{SH. Chen, S. Ponnusamy,
X. Wang }{On properties of solutions of the $p$-harmonic
equation}

\section{Introduction and Preliminaries}\label{csw-sec1}

A complex-valued function $f=u+iv$ in a simply connected domain $\Omega \subseteq\mathbb{C}$ is called
\textit{p-harmonic} if $u$ and $v$ are $p$-harmonic in $\Omega$, i.e.
$f$ satisfies the $p$-harmonic equation $\Delta ^pf=0$, where
$$\Delta ^pf= \underbrace{\Delta\cdots\Delta}_{p}f,
$$
where $p$ is a positive integer and $\Delta$ represents the Laplacian
operator
$$\Delta:=4\frac{\partial^{2}}{\partial
z\partial\overline{z}}=\frac{\partial ^2}{\partial x^2} +
\frac{\partial ^2}{\partial y^2}.
$$
Throughout this paper we consider $p$-harmonic mappings of the unit disk
$\ID=\{z\in \mathbb{C}:\, |z|<1\}.$ Obviously, when $p=1$ (resp. $p=2$), $f$ is harmonic (resp.
biharmonic). The properties of harmonic \cite{Clunie-Small-84,D} and biharmonic \cite{AA,A,AM,Kh,La}
mappings have been investigated by many authors.
Concerning $p$-harmonic mappings, we easily have the following characterization.

\begin{prop}\label{prop1a}
A mapping $f $ is $p$-harmonic in $\ID$ if and only if $f$ has the following representation:
\be\label{eq2.0x}
f(z)=\sum_{k=1}^{p}|z|^{2(k-1)}G_{p-k+1}(z),
\ee
where $G_{p-k+1}$ is harmonic for each $k\in \{1,\ldots ,p\}$.
\end{prop}
\bpf
We only need to prove the necessity since the proof for the
sufficiency part is obvious. Again, as the cases $p=1,2$ are
well-known, it suffices to prove the result for $p\geq 3$. We shall
prove the proposition by the method of induction. So, we assume that
the proposition is true for $p=n \,(\geq 3).$

Let $F$ be an $(n+1)$-harmonic mapping in $\ID$.
By assumption, $\Delta F$ is $n$-harmonic and so can be represented as
$$
\Delta F(z)=\sum_{k=1}^{n}|z|^{2(k-1)}G_{n-k+1}(z),
$$
where $G_{n-k+1}$ ($1\leq k\leq n$) are harmonic functions with
$$G_{n-k+1}(z)=a_{0,n-k+1}+\sum_{j=1}^{\infty}a_{j,n-k+1}z^{j}+
\sum_{j=1}^{\infty}\overline{b}_{j,n-k+1}\overline{z}^{j}
\quad \mbox{for $k\in \{1, \ldots , n\}$.}
$$
Then
\begin{eqnarray*}
\int_{0}^{z}\int_{0}^{\overline{z}}\Delta
F\,d\overline{z}\,dz&=&\sum_{k=1}^{n}|z|^{2k}T_{p-k+1}(z)+g(z),
\end{eqnarray*}
where
$$T_{p-k+1}(z)=\sum_{k=1}^{n}\left(\frac{a_{0,n-k+1}}{k^{2}}+
\sum_{j=1}^{\infty}\frac{a_{j,n-k+1}}{k(k+j)}z^{j}+
\sum_{j=1}^{\infty}\frac{\overline{b}_{j,n-k+1}}{k(k+j)}\overline{z}^{j}\right)
$$
and $ g$ is a harmonic function in $\ID$. A rearrangement of the
series in the sum shows that (\ref{eq2.0x}) holds for $p=n+1$.
\epf

We remark that the representation (\ref{eq2.0x}) continues to hold even if $f$ is $p$-harmonic in
a simply connected domain $\Omega$.


For a sense-preserving $C^1$-mapping (i.e. continuously differentiable), we let
$$\lambda_{f}=|f_{z}|-|f_{\overline{z}}|  ~\mbox{ and }~
\Lambda_{f}=|f_{z}|+|f_{\overline{z}}|
$$
so that the Jacobian $J_{f}$ of $f$ takes the form
$$J_{f}=\lambda_{f}\Lambda_{f}=|f_{z}|^{2}-|f_{\overline{z}}|^{2}>0.
$$
In \cite{Al-Amri-Mocanu80}, the authors obtained sufficient
conditions for the univalence of $C^{1}$-functions. Now we introduce
the concepts of starlikeness and convexity of $C^{1}$-functions.

\bdefe
A $C^1$-mapping $f$ with $f(0)=0$ is called starlike if $f$ maps
$\ID$ univalently onto a domain $\Omega$ that is starlike with
respect to the origin, i.e. for every $w\in \Omega$ the line segment
$[0,w]$ joining $0$ and $w$ is contained in $\Omega$. It is known
that $f$ is starlike if it is sense-preserving, $f(0)=0$, $f(z)\neq
0$ for all $z\in \ID\setminus \{0\}$ and
$$\frac{\partial}{\partial t}\big(\arg f(re^{it})\big):=
{\rm Re} \left ( \frac{Df(z)}{f(z)}\right ) >0
\quad \mbox{for all $z=re^{it}\in \ID\setminus \{0\}$},
$$
where $Df=zf_{z}-\overline zf_{\overline z}$ (cf. \cite[Theorem
1]{Mocanu80}).
\edefe



\bdefe \label{def2}
Let  $f$ and $Df$ belong to $C^1(\ID)$. Then we say that $f$ is {\it convex} in $\ID$
if it is sense-preserving, $f(0)=0$,
$f(z)\cdot Df(z)\neq0$ for all $z\in \ID\setminus \{0\}$ and
$$
{\rm Re}\left (\frac{D^2f(z)}{Df(z)}\right )>0 \quad \mbox{for all $z\in \ID\setminus \{0\}$}.
$$
\edefe

As $\arg Df(re^{it})$ represents the argument of the outer normal to the curve
$C_r=\{f(re^{i\theta}):\, 0\leq \theta <2\pi\}$ at the point $f(re^{it})$, the last condition gives that
$$\frac{\partial}{\partial t}\big(\arg Df(re^{it})\big)
= {\rm Re}\left (\frac{D^2f(z)}{Df(z)}\right )>0 \quad \mbox{for all $z=re^{it}\in \ID\setminus \{0\}$},
$$
showing that the curve $C_r$ is convex for each $r\in(0, 1)$ (see \cite[Theorem 2]{Mocanu80}).
Non-analytic starlike and convex  functions were studied by
Mocanu in \cite{Mocanu80}. Harmonic starlike and harmonic convex functions were systematically
studied by  Clunie and  Sheil-Small \cite{Clunie-Small-84}, and these two classes of
functions have been studied extensively by many authors. See for instance, the
book by Duren \cite{D} and the references therein.



The complex differential operator
$$D=z\frac{\partial}{\partial z}-\overline{z}\frac{\partial}{\partial \overline{z}}
$$
defined by Mocanu \cite{Mocanu80} on the class of complex-valued $C^1$-functions
satisfies the usual product rule:
$$D(af+bg)=aD(f)+bD(g) ~\mbox{ and }~D(fg)=fD(g)+gD(f),
$$
where $a,b$ are complex constants, $f$ and $g$ are $C^1$-functions. The operator $D$
possesses a number of interesting properties. For instance, the operator $D$
preserves both  harmonicity and  biharmonicity (see also \cite{AM}).
In the case of $p$-harmonic mappings, we also have the following property of
the operator $D$.

\begin{prop}\label{prop1}
$D$ preserves $p$-harmonicity.
\end{prop}
\bpf Let $f$ be a $p$-harmonic mapping  with the form
$$f(z)=\sum_{k=1}^{p}|z|^{2(k-1)}G_{p-k+1}(z),
$$
where each $G_{p-k+1}(z)$ is harmonic in $\mathbb{D}$ for $k\in \{1,
\ldots, p\}$. As $D\big (|z|^2\big )=0$, the product rule shows that
$D\big (|z|^{2(k-1)}\big )=0$ for each $k\in \{1, \ldots, p\}$. In
view of this and the fact that $D$ preserves harmonicity gives that
\begin{eqnarray*}
 D(f(z))&=&\sum_{k=1}^{p}\Big
[|z|^{2(k-1)}D(G_{p-k+1}(z))+ D
(|z|^{2(k-1)})G_{p-k+1}(z)\Big ]\\
&=& \sum_{k=1}^{p}
|z|^{2(k-1)}D(G_{p-k+1}(z)).
\end{eqnarray*}
 \epf

%
%
%

One of the aims of this paper is to generalize the main results of
Abdulhadi, et. al.  \cite{AM} to the case of
$p$-harmonic mappings. The corresponding generalizations are Theorems \ref{thm1} and \ref{thm2}.


The classical theorem of Landau for bounded analytic functions states that
if $f$ is analytic in $\mathbb{D}$ with $f(0)=f'(0)-1=0$, and $|f(z)|<M$ for $z\in \ID$,
then $f$ is univalent in the disk $\mathbb{D}_\rho :=\{z\in \mathbb{C}:\;
|z|<\rho \} $ and in addition, the range $f(\mathbb{D}_\rho )$
contains a disk of radius $M\rho ^2$ (cf. \cite{L}), where
$$\rho =\frac{1}{M+\sqrt{M^2 -1}}.
$$
Recently, many authors considered Landau's theorem for planar
harmonic mappings (see for example,
\cite{HG,CPW,CPW1,DN-04,Armen-06,Liu,Hu}) and biharmonic mappings
(see \cite{AA,CPW0,CPW,Liu0}). In Section \ref{sec-5}, we consider
Landau's theorem for $p$-harmonic mappings with the form $D(f)$ when
$f$ belongs to certain classes of $p$-harmonic mappings. Our results
are Theorems \ref{thm3.x} and \ref{thm4}.

In a series of papers the second author with Yanagihara and
Vasudevarao (see \cite{S,S1,Y1,Y2}) have discussed the regions of
variability for certain classes of univalent analytic functions in
 $\ID$. In Section \ref{sec-4} (see Theorem
\ref{thm3.1}), we solve a related problem for certain $p$-harmonic
mappings. Finally, in Section \ref{sec6}, we present explicit upper estimates for Bloch norm for
bi- and tri-harmonic mappings (see Corollaries  \ref{cor1a} and \ref{cor2a}).

\section{Lemmas}\label{sec-2}

For the proofs of our main results we require a number of lemmas. We begin to recall
the following version of Schwarz lemma due to Heinz (\cite[Lemma]{He}) and
Colonna \cite[Theorem 3]{Co}, see also \cite{HG,CPW,CPW1}.

\begin{Lem}\label{LemA} Let $f$ be a harmonic mapping of $\mathbb{D}$
such that $f(0)=0$ and $f(\mathbb{D})\subset\mathbb{D}$. Then
$$
|f(z)|\leq\frac{4}{\pi}\arctan|z|\leq\frac{4}{\pi}|z|\;\;
\mbox{for}\;\; z\in\mathbb{D}
$$
and
$$
 \Lambda_{f}(z)\leq\frac{4}{\pi}\frac{1}
{(1-|z|^{2})}\;\; \mbox{for}\;\; z\in\mathbb{D}.
$$
\end{Lem}


\begin{Lem}{\rm (\cite[Lemma 2.1]{Liu})}\label{LemB}
Suppose that $f(z)=h(z)+\overline{g(z)}$ is a harmonic mapping of
$\mathbb{D}$ with $h(z)=\sum_{n=1}^{\infty}a_{n}z^{n}$ and
$g(z)=\sum_{n=1}^{\infty}b_{n}z^{n}$ for $z\in\mathbb{D}.$
If $J_{f}(0)=1$ and $|f(z)|<M,$ then
$$
|a_{n}|,\ |b_{n}|\leq\sqrt{M^{2}-1},\ n=2, 3, \ldots,
$$
$$
|a_{n}|+ |b_{n}|\leq\sqrt{2M^{2}-2},\ n=2, 3, \ldots
$$
and
\be \label{eqe}
\lambda_{f}(0)
 \geq \lambda_{0}(M) := \begin{cases}
\ds \frac{\sqrt{2}}{\sqrt{M^{2}-1}+\sqrt{M^{2}+1}}&
\mbox{ if } \ds 1\leq M\leq M_0, \\[4mm]
\ds \frac{\pi}{4M} &  \mbox{ if }
\ds M> M_0,\end{cases}
\ee where $M_0=\frac{\pi}{2\sqrt[4]{2\pi^{2}-16}} \approx 1.1296$.
\end{Lem}

The following lemma concerning coefficient estimates for harmonic mappings is crucial
in the proofs of Theorems \ref{thm1} and \ref{thm2}.
This lemma has been proved by the authors in \cite{ CPW2} with an additional assumption
that $f(0)=0$. However, for the sake of clarity, we present a slightly different proof
than that in \cite{ CPW2}.

\begin{Lem}\label{Lem1}
Let $f=h+\overline{g}$ be a harmonic mapping of $\mathbb{D}$ such that  $|f(z)|<
M$ with $h(z)=\sum_{n=0}^{\infty}a_{n}z^{n}$ and $g(z)=\sum_{n=1}^{\infty}b_{n}z^{n}$.
Then $|a_{0}|\leq M$ and for any $n\geq 1$
\be\label{eq2.4}
|a_{n}|+|b_{n}|\leq \frac{4M}{\pi}.
\ee
The estimate $(\ref{eq2.4})$ is sharp. The extremal functions are $f(z)\equiv M$ or
$$f_{n}(z)=\frac{2M\alpha }{\pi}\arg \left ( \frac{1+\beta z^{n}}{1-\beta z^{n}}\right),
$$
where $|\alpha|=|\beta |=1$.
\end{Lem}
\bpf
Without loss of generality, we assume that $|f(z)|< 1$. For $\theta\in [0, 2\pi )$, let
$$v_{\theta}(z)={\rm Im\,}(e^{i\theta}f(z))
$$
and observe that
$$v_{\theta}(z)  = {\rm Im\,}( e^{i\theta}h(z)+\overline{e^{-i\theta}g(z)})
={\rm Im\,}(e^{i\theta}h(z)-e^{-i\theta}g(z)).
$$
Because $|v_{\theta}(z)| <1$, it follows that
$$e^{i\theta}h(z)-e^{-i\theta}g(z)\prec K(z)
=\lambda+\frac{2}{\pi}\log \left (\frac{1+z\xi}{1-z} \right ),
$$
where $\xi=e^{-i\pi\mbox{Im}(\lambda)}$ and
$\lambda=e^{i\theta}h(0)-e^{-i\theta}g(0)$. The superordinate
function $K(z)$ maps $\ID$ onto a convex domain with $K(0)=\lambda$
and $K'(0)=\frac{2}{\pi}(1+\xi)$, and therefore, by a theorem of
Rogosinski \cite[Theorem 2.3]{Rago-43} (see also \cite[Theorem
6.4]{Du}), it follows that
$$|a_{n}-e^{-2i\theta}b_{n}|\leq\frac{2}{\pi}|1+\xi|\leq \frac{4}{\pi} \quad \mbox{ for  $n=1, 2,\ldots $}
$$
and the desired inequality (\ref{eq2.4}), with $M=1$, is a consequence of the
arbitrariness of $\theta$ in $[0, 2\pi )$.

For the proof of sharpness part, consider the functions
$$f_{n}(z)=\frac{2M\alpha}{\pi}\mbox{Im}\left(\log\frac{1+\beta z^{n}}{1-\beta z^{n}}\right),
\quad |\alpha|=|\beta |=1,
$$
whose values are confined to a diametral segment of the disk $\ID_M 
$. Also,
$$f_{n}(z)=\frac{2M\alpha }{i\pi}\left(\sum_{k=1}^{\infty}\frac{1}{2k-1}(\beta z^{n})^{2k-1}-
\sum_{k=1}^{\infty}\frac{1}{2k-1}(\overline{\beta}
\overline{z}^{n})^{2k-1}\right),
$$
which gives
$$|a_{n}|+|b_{n}|=\frac{4M}{\pi}.
$$
The proof of the lemma is complete.
\epf


As an immediate consequence of Lemmas \Ref{LemB} and \Ref{Lem1}, we have

\begin{cor}\label{cor2}
Let $f=h+\overline{g}$ be a harmonic mapping of $\mathbb{D}$
with $h(z)=\sum_{n=1}^{\infty}a_{n}z^{n}$,
$g(z)=\sum_{n=1}^{\infty}b_{n}z^{n}$ and $|f(z)|\leq M$. If
$J_{f}(0)=1$ and $M\geq\frac{\pi}{\sqrt{\pi^{2}-8}}$, then for any
$n\geq 2$,
$$
|a_{n}|+ |b_{n}|\leq \frac{4M}{\pi}\leq\sqrt{2M^{2}-2}.
$$
\end{cor}

\section{The convexity and the starlikeness}\label{sec-3}


The following simple result can be used to generate (harmonic) starlike and convex functions.

\begin{thm}\label{thm1}
Let $f$ be a univalent $p$-harmonic mapping with the form
$$f(z)=G(z) \sum_{k=1}^{p}\lambda_{k}|z|^{2(k-1)},
$$
where $G$ is a locally univalent harmonic mapping and $\lambda_{k}\ (k=1, \ldots, p)$ are complex constants.
Then we have the following:
\begin{enumerate}
\item [{\rm \textbf{(a)}}] $\ds \frac{D(f)}{f}=\frac{D(G)}{G}$ and $\ds \frac{D^2(f))}{D(f)}=\frac{D^2(G))}{D(G)}.
$
\item [{\rm \textbf{(b)}}]  $f$ is convex $($resp. starlike$)$ if and
only if $G$ is convex $($resp. starlike$)$.
\end{enumerate}
\end{thm}
\bpf \textbf{(a)} The two equalities are immediate consequences of
the formula
$$D\Big (G(z)\sum_{k=1}^{p}\lambda_{k}|z|^{2(k-1)}\Big )=D(G(z))\sum_{k=1}^{p}\lambda_{k}|z|^{2(k-1)}.
$$
So, we omit the details.

\textbf{(b)} It suffices to prove the case of convexity since the proof for the starlikeness is similar.

Let $z=re^{it}$, where $0<r<1$ and $0\leq t<2\pi$. Then
$$f(z)=G(z) \sum_{k=1}^{p}\lambda_{k}|z|^{2(k-1)}=G(re^{i\theta}) \sum_{k=1}^{p}\lambda_{k}r^{2(k-1)},
$$
so that
$$\frac{\partial f(re^{it})}{\partial t}=\frac{\partial G(re^{it})}{\partial t} \,
\sum_{k=1}^{p}\lambda_{k}r^{2(k-1)}
$$
and
$$\frac{\partial^{2}f(re^{it})}{\partial t^{2}}= \frac{\partial^{2}G(re^{it})}{\partial t^{2}}\,
\sum_{k=1}^{p}\lambda_{k}r^{2(k-1)}.
$$
Therefore Part \textbf{(a)} yields
\begin{align*}
\frac{\partial}{\partial t} \Big (\arg\frac{\partial f(re^{it})
}{\partial t}\Big  ) ={\rm  Re\,} \Big (\frac{D^2(f)}{D(f)}\Big  )
={\rm  Re\,} \Big (\frac{D^2(G)}{D(G)}\Big  ) =\frac{\partial
}{\partial t} \Big (\arg\frac{\partial G(re^{it})}{\partial t}\Big ),
\end{align*}
from which the proof of Part \textbf{(b)} of this theorem follows.
\epf

As an immediate consequence of Theorem \ref{thm1}{\bf (a)}, we easily have the
following.

\begin{cor}\label{cor1}
Let $f$ be a univalent  $p$-harmonic mapping defined as in Theorem {\rm \ref{thm1}}.
 If $f$ is convex and $D(f)$ is univalent, then $D(f)$ is
starlike.
\end{cor}

Abdulhadi, et. al. \cite[Theorem 1]{AM} discussed the
univalence and the starlikeness of biharmonic mappings in $\mathbb{D}$.
A natural question is whether \cite[Theorem 1]{AM} holds for $p$-harmonic mappings.
The following result gives a partial answer to this problem.

\begin{thm}\label{thm2}
Let $f$ be a $p$-harmonic mapping of $\mathbb{D}$ satisfying
$f(z)=|z|^{2(p-1)}G(z)$, where $G$ is  harmonic, orientation
preserving and starlike.  Then $f$ is starlike univalent.
\end{thm}
\bpf We see that the Jacobian $J_{f}$ of $f$ is
\begin{align*}
J_{f}&=|f_{z}|^{2}-|f_{\overline{z}}|^{2}\\
&=|z|^{4(p-1)}(|G_{z}|^{2}-|G_{\overline{z}}|^{2})+2(p-1)|z|^{4p-6}|G|^{2}\mbox{Re\,}\left (\frac{D(G)}{G}\right  )\\
&\geq |z|^{4(p-1)}(|G_{z}|^{2}-|G_{\overline{z}}|^{2}).
\end{align*}
Hence $J_{f}(z)>0$ when $0<|z|<1$ and obviously, $J_{f}(0)=0.$ The
univalence of $f$ follows from a standard argument as in the proof of \cite[Theorem 1]{AM}.
Finally, Theorem \ref{thm1} implies that $f$ is starlike.
\epf

\section{The Landau theorem}\label{sec-5}

We now discuss the existence of the Laudau constant for two
classes of $p$-harmonic mappings.

\begin{thm}\label{thm3.x}
Let $f(z)=\sum_{k=1}^{p}|z|^{2(k-1)}G_{p-k+1}(z)$ be a $p$-harmonic
mapping of $\mathbb{D}$ satisfying $\Delta
G_{p-k+1}(z)=f(0)=G_{p}(0)=J_{f}(0)-1=0$ and for any
$z\in\mathbb{D}$, $|G_{p-k+1}(z)|\leq M$, where $M\geq 1$. Then
there is a constant $\rho$ $(0<\rho<1)$ such that $D(f)$ is
univalent in $\mathbb{D}_{\rho}$, where $\rho$ satisfies the
following equation:
$$\lambda_{0}(M)-\frac{T(M)}{(1-\rho)^{2}}\sum_{k=2}^{p}(2k-1)\rho^{2(k-1)}
-\sum_{k=1}^{p}\frac{2T(M)\rho^{2k-1}}{(1-\rho)^3}
-\frac{16M}{\pi^{2}}s_{0}\arctan\rho=0
$$
with
$$s_{0}= \left ( \frac{\sqrt{17}-1}{\sqrt{17}-3} \right ) \sqrt{\frac{2}{5-\sqrt{17}}}~
\approx 4.1996,
$$
\be\label{eqtm1}
T(M)=\begin{cases}
\sqrt{2M^{2}-2}  & \mbox{ if } \ds 1\leq M\leq
M_1:=\frac{\pi}{\sqrt{\pi^{2}-8}}\approx 2.2976\\
\ds \frac{4M}{\pi} & \mbox{ if }\ds
M>M_1
\end{cases}
\ee
and $\lambda_{0}(M)$ is given by $(\ref{eqe})$.  Moreover, the range
$D(f)(\mathbb{D}_{\rho})$ contains a univalent disk
$\mathbb{D}_{R}$, where
$$R=\rho\Big[\lambda_{0}(M)-\sum_{k=2}^{p}\frac{T(M)\rho
^{2(k-1)}}{(1-\rho )^{2}}- \frac{16M}{\pi^{2}}s_{0}\arctan\rho\Big]
.
$$
\end{thm}
\bpf
For each $k\in\{1, 2, \ldots, p\}$, let
$$G_{p-k+1}(z)=a_{0,p-k+1}+\sum_{j=1}^{\infty}a_{j,p-k+1}z^{j}
+\sum_{j=1}^{\infty}\overline{b}_{j,p-k+1}\overline{z}^{j},
$$
where   $a_{0,p}=0$. We define the function $H$ as
$$H=D\left  (\sum_{k=1}^{p}|z|^{2(k-1)}G_{p-k+1}\right )
=\sum_{k=1}^{p}|z|^{2(k-1)}D(G_{p-k+1}).
$$
Using Lemmas \Ref{LemB}, \Ref{Lem1} and Corollary \ref{cor2}, we
have
$$|a_{n,p}|+|b_{n,p}|\leq T(M),
$$
where $T(M)$ is given by (\ref{eqtm1}), and
$$|a_{j,p-k+1}|+|b_{j,p-k+1}|\leq\frac{4M}{\pi}
$$
for $j\geq 1,$ $n\geq 2$ and $2\leq k\leq p.$

We observe that
$$J_{f}(0)=|(G_{p})_{z}(0)|^{2}-|(G_{p})_{\overline{z}}(0)|^{2}=J_{G_{p}}(0)=1
$$
and hence by Lemmas \Ref{LemA} and \Ref{LemB}, we have
$$\lambda_{f}(0)\geq \lambda_{0}(M),
$$
where $\lambda_{0}(M)$ is given by $(\ref{eqe})$. Now, we define
$$q(x)=\frac{2-x^{2}}{(1-x^{2})x}\;\,(0<x<1).
$$
Then there is an $r_0=\sqrt{\frac{5-\sqrt{17}}{2}}\approx0.66$ such
that
$$q(r_{0})=\min_{0<x<1}q(x)=\left (\frac{\sqrt{17}-1}{\sqrt{17}-3}\right )\sqrt{\frac{2}{5-\sqrt{17}}}
~=s_0.
$$
For each $\theta\in[0,2\pi)$, the function
$$G_{\theta}(z)=(G_{p})_{z}(z)-(G_{p})(0)+((G_{p})_{\overline{z}}(z)-(G_{p})_{\overline{z}}(0))e^{i(\pi-2\theta)}\;\,
$$
is clearly a harmonic mapping of $\mathbb{D}$ and satisfies
$G_{\theta}(0)=0.$ Moreover, it follows from Lemma \Ref{LemA} that
$$
\Lambda_{G_{p}}(z)\leq\frac{4M}{\pi}\frac{1}{1-|z|^{2}} ~\mbox{ for
$z\in\mathbb{D}$.}
$$
In particular, this observation yields that \be\label{eq(2.4)}
|G_{\theta}(z)|\leq\Lambda_{G_{p}}(z)+\Lambda_{G_{p}}(0)\leq
\frac{4M}{\pi}\Big (1+\frac{1}{1-|z|^{2}}\Big
)=\frac{4M}{\pi}|z|q(|z|)
\ee for all $z\in\mathbb{D}$.

Since $xq(x)-1=\frac{1}{1-x^{2}}$ is an increasing function
in the interval $(0,1)$, the inequality (\ref{eq(2.4)}) shows that for any $z\in
\mathbb{D}_{r_0}$,
$$|G_\theta(z)|\leq \frac{4M}{\pi}m_0,
$$
where $m_0=(2- r_{0}^{2})/(1-r_{0}^{2}). $
Next, we consider the mapping $F$ defined on $\mathbb{D}$ by
$$F(z)=\frac{\pi}{4Mm_0}G_{\theta}(r_0z).
$$
Applying Lemma \Ref{LemA} to the function $F(z)$ yields
that for $z\in \mathbb{D}_{r_0}$,
$$|{G}_{\theta}(z)|\leq\frac{16M}{\pi^{2}}m_0\arctan \left (\frac{|z|}{r_0}\right )
\leq \frac{16M}{\pi^{2}}s_{0}\arctan |z|,
$$
where $s_{0}= m_0/r_0$.

Now, we fix $\rho$ with $\rho\in (0,1)$. To prove the univalency of $H$,
we choose two distinct points $z_{1}, z_{2}$ in $\ID_
\rho$. Let $\gamma =\{(z_{2}-z_{1})t+z_{1}:\;0\leq t\leq 1\}$ and
$z_{2}-z_{1}=|z_{1}-z_{2}|e^{i\theta}$.
We find that

\vspace{6pt}

$|H(z_1)-H(z_2)|$
\begin{eqnarray*}
 &=& \Big |\int_{\gamma }H_z(z)\,dz+H_{\overline{z}}(z)\,d\overline{z}\Big |\\
&\geq & \Big |\int_{\gamma}(G_{p})_{z}(0)\,dz-(G_{p})_{\overline{z}}(0)\,d\overline{z}\Big  |\\
&&-\Big |\int_{\gamma}\sum_{k=2}^{p}|z|^{2(k-1)}[z(G_{p-k+1})_{z^{2}}(z)\,dz
-\overline{z}(G_{p-k+1})_{\overline{z}^{2}}(z)\,d\overline{z}]\Big  |\\
&& -\Big |\int_{\gamma}\sum_{k=2}^{p}(k-1)|z|^{2(k-2)}[z^{2}(G_{p-k+1})_{z}(z)\,d\overline{z}
-\overline{z}^{2}(G_{p-k+1})_{\overline{z}}(z)\,dz]\Big  |\\
&& -\Big |\int_{\gamma}\sum_{k=2}^{p}k|z|^{2(k-1)}[(G_{p-k+1})_{z}(z)\,dz-(G_{p-k+1})_{\overline{z}}(z)
\,d\overline{z}]\Big  |\\
&& -\Big|\int_{\gamma}[(G_{p})_{z}(z)-(G_{p})_{z}(0)]\,dz-[(G_{p})_{\overline{z}}(z)-
(G_{p})_{\overline{z}}(0)]\,d\overline{z}\Big |\\
&\geq& |z_{1}-z_{2}|\Big \{\lambda_{f}(0)-|{G}_{\theta}(\rho)| \Big .\\
&& -\sum_{k=1}^{p}\rho^{2(k-1)}\sum_{n=2}^\infty
n(n-1)(|a_{n,p-k+1}|+|b_{n,p-k+1}|)\rho^{n-1}\\
&&-\Big .\sum_{k=2}^{p}(2k-1)\rho^{2(k-2)}\sum_{n=1}^{\infty}n(|a_{n,p-k+1}|+|b_{n,p-k+1}|)\rho^{n+1} \Big  \}\\
&>& |z_{1}-z_{2}|\Big
[\lambda_{0}(M)-\frac{T(M)}{(1-\rho)^{2}}\sum_{k=2}^{p}(2k-1)\rho^{2(k-1)} \Big.
\\
&& \Big .-\sum_{k=1}^{p}\frac{2T(M)\rho^{2k-1}}{(1-\rho)^3}
-\frac{16M}{\pi^{2}}s_{0}\arctan\rho \Big  ].
\end{eqnarray*}

Let
$$P(\rho)=\lambda_{0}(M)-\frac{T(M)}{(1-\rho)^{2}}\sum_{k=2}^{p}(2k-1)\rho^{2(k-1)}
-\sum_{k=1}^{p}\frac{2T(M)\rho^{2k-1}}{(1-\rho)^3}
-\frac{16M}{\pi^{2}}s_{0}\arctan\rho.
$$
Then it is easy to verify that $P(\rho)$ is a decreasing function on the interval $(0,1)$,
$$\lim_{\rho\rightarrow 0+}P(\rho)=\lambda_{0}(M) ~\mbox{ and }~ \lim_{\rho\rightarrow 1-}
P(\rho)=-\infty.
$$
Hence there exists a unique $\rho_{0}$ in $(0,1)$ satisfying $P(\rho _0)=0.$
This observation shows that $|H(z_1)-H(z_2)|>0$ for
arbitrary two distinct points $z_{1}, z_{2}$ in $|z|<\rho_{0}$
which proves the univalency of $D(F)$ in $\ID_{\rho_{0}}$.

For any $z$ with $|z|=\rho_{0} $, we have
\begin{eqnarray*}
|H(z)| &= & \Big|\sum_{k=1}^{p}|z|^{2(k-1)}[z(G_{p-k+1})_{z}(z)-\overline{z}(G_{p-k+1})_{\overline{z}}(z)] \Big  | \\
&\geq& \Big| z(G_{p})_{z}(0)-\overline{z}(G_{p})_{\overline{z}}(0)\Big  |\\
&&- \Big |z[(G_{p})_{z}(z)-(G_{p})_{z}(0)]
-\overline{z}[(G_{p})_{\overline{z}}(z)- (G_{p})_{\overline{z}}(0)]\Big  |\\
&& -\Big|\,\sum_{k=2}^{p}|z|^{2(k-1)}[z(G_{p-k+1})_{z}(z)-\overline{z}(G_{p-k+1})_{\overline{z}}(z)]\Big  |\\
& \geq&\rho_{0} \Big [\lambda_{0}(M)-\sum_{k=2}^{p}\frac{T(M)\rho_{0}
^{2(k-1)}}{(1-\rho_{0} )^{2}}- \frac{16M}{\pi^{2}}s_{0}\arctan\rho_{0}\Big ]\\
& =& R
\end{eqnarray*}
and the proof of the theorem is complete. \epf

\begin{table}
\center{
\begin{tabular}{|c|c||c|c|c|c|}
\hline
$M$ & $p$ &$\rho=\rho(M,p)$ & $R=R(M,\rho (M,p))$  & $\rho '$ & $R'$
\\
\hline
1.1296& 2&0.0714741  & 0.0101601 &0.0420157 & 0.00945379
\\
2& 2&0.0206783  & 0.00227639 &0.0139439 & 0.00164502
\\
2.2976& 2&0.0155966&0.00151523 &0.0106132  & 0.00108021
\\
3& 2&0.00922255  &0.00067425  &0.00626141 &0.000482413
\\
1.1296& 3&0.071463    & 0.0101647 & -- & --
\\
2& 3&0.0206782  &0.00227641 & -- & --
\\
2.2976& 3&0.0155966  &0.00151523 & -- & --
\\
3&3 &0.00922254 &0.000674251  & -- & --
\\
1.1296& 4&0.0714629    &0.0101647  & -- & --
\\
2& 4&0.0206782  &0.00227641 & -- & --
\\
2.2976& 4&0.0155966  &0.00151523 & -- & --
\\
3&4 &0.00922254 &0.000674251  & -- & --
\\
\hline
\end{tabular}
}
\bigskip
\caption{Values of $\rho$ and $R$ for Theorem \ref{thm3.x} for $p=2$,
and the corresponding values of $\rho '$ and $R '$ of \cite[Theorem 1.1]{CPW0} (for $p=2$)
\label{table1}
}
\end{table}

From Table \ref{table1}, we see that Theorem \ref{thm3.x} improves Theorem~1.1 of \cite{CPW0} for the case $p=2$,
and the results for the rest of the values of $p$ are new. In Table \ref{table1}, third and fourth columns
refer to values obtained from Theorem \ref{thm3.x} for cases $p=2,3,4$ for certain choices of
$M$, while the right two columns correspond to the values obtained from \cite[Theorem 1.1]{CPW0} for the case
$p=2$.

\begin{thm}\label{thm4}
Let $f(z)=|z|^{2(p-1)}G(z)$ be a $p$-harmonic mapping of
$\mathbb{D}$ satisfying $G(0)=J_{G}(0)-1=0$ and $|G(z)|\leq M$, where
$M\geq 1$ and $G$ is harmonic.  Then there is a constant $\rho$
$(0<\rho<1)$ such that $D(f)$ is univalent in $\mathbb{D}_{\rho}$,
where $\rho$ satisfies the following equation:
$$\lambda_{0}(M)-\frac{48M}{\pi^{2}}s_{0}\arctan\rho-\frac{2T(M)\rho}{(1-\rho)^{3}}=0,
$$
where the constants $s_{0}$, $\lambda_{0}(M)$ and $T(M)$ are the same as in
Theorem {\rm \ref{thm3.x}}. Moreover, the range $D(f)(\mathbb{D}_{\rho})$ contains a univalent
disk $\mathbb{D}_{R}$, where
$$R =\rho^{2p-1}\Big [\lambda_{0}(M)-\frac{16M}{\pi^{2}}s_{0}\arctan\rho
\Big ].
$$
Especially, if $M=1$, then $G(z)=z$, i.e. $f(z)=|z|^{2(p-1)}z$ which
is univalent in $\mathbb{D}$.
\end{thm}
\bpf Let $G(z)=\sum_{n=1}^{\infty}a_{n}z^{n}+\sum_{n=1}^{\infty}\overline{b}_{n}\overline{z}_{n}.$
 Using Lemmas \Ref{LemB}, \Ref{Lem1} and Corollary \ref{cor2}, we have
$$|a_{n}|+|b_{n}|\leq T(M) ~\mbox{ for  $n\geq 2$.}~
$$
Note that
$$J_{G}(0)=|a_{1}|^{2}-|b_{1}|^{2}=1
$$
and hence, by Lemmas \Ref{LemA} and \Ref{LemB}, we have
$$\lambda_{G}(0)\geq \lambda_{0}(M).
$$
Next, we set $H=D(f)=|z|^{2(p-1)}D(G)$ and  fix $\rho$ with
$\rho\in (0,1)$. To prove the univalency of $f$, we choose two distinct
points $z_{1}, z_{2}$ in $\ID_ \rho$. Let $\gamma
=\{(z_{2}-z_{1})t+z_{1}:\;0\leq t\leq 1\}$ and
$z_{2}-z_{1}=|z_{1}-z_{2}|e^{i\theta}$.
Then
\begin{eqnarray*}
|H(z_{1})-H(z_{2})|&=&
\left  |\int_{[z_{1},z_{2}]}H_{z}(z)\,dz+H_{\overline{z}}(z)\,d\overline{z}\right |\\
&=&\left |\int_{[z_{1},z_{2}]}p|z|^{2(p-1)}(G_{z}(z)\,dz-G_{\overline{z}}(z)\,d\overline{z})\right .\\
&& +|z|^{2(p-1)}(zG_{z^{2}}(z)\,dz-\overline{z}G_{\overline{z}^{2} }(z)\,d\overline{z}) \\
&& \left .+ (p-1)|z|^{2(p-2)}(z^{2}G_{z}(z)\,d\overline{z}-\overline{z}^{2}G_{\overline{z}}(z)\,dz) \right  |\\
&\geq& \left |\int_{[z_{1},z_{2}]}\left [G_{z}(0)(p|z|^{2(p-1)}\,dz+(p-1)|z|^{2(p-2)}z^{2}\,d\overline{z})
\right .\right .\\
&& \left .\left . - G_{\overline{z}}(0)(p|z|^{2(p-1)}\,d\overline{z} -(p-1)|z|^{2(p-2)}\overline{z}^{2}\,dz)\right  ]\frac{}{} \right  |\\
& & -p\left |\int_{[z_{1},z_{2}]}|z|^{2(p-1)}\big [(G_{z}(z)-G_{z}(0))\,dz-(G_{\overline{z}}(z)
-G_{\overline{z}}(0)) \,d\overline{z}]\right  |\\
&&- \left|(p-1)\int_{ [z_{1},z_{2}]}|z|^{2(p-1)}\Big[\frac{z}{\overline{z}}(G_{z}(z)-G_{z}(0))\,d\overline{z}\right .\\
&&\left .-\frac{\overline{z}}{z}(G_{\overline{z}}(z)-G_{\overline{z}}(0))\,dz\Big ]
\right |\\
&&-\left |\int_{[z_{1},z_{2}]}|z|^{2(p-1)}(zG_{z^{2}}(z)\,dz-\overline{z}G_{\overline{z}^{2}}(z)\,d\overline{z})
\right |\\
&\geq& |z_{1}-z_{2}|\Big (\int_{0}^{1}|z|^{2(p-1)}dt\Big )\left \{\lambda_{0}(M)- \frac{48M}{\pi^{2}}s_{0}\arctan\rho
\right .\\
&& \left . -\sum_{n=2}^{\infty}n(n-1)(|a_{n}|+|b_{n}|)\rho^{n-1}\right \}\\
&>& |z_{1}-z_{2}|\Big (\int_{0}^{1}|z|^{2(p-1)}dt\Big )\left [\lambda_{0}(M)-
\frac{48M}{\pi^{2}}s_{0}\arctan\rho-\frac{2T(M)\rho}{(1-\rho)^{3}}\right ].
\end{eqnarray*}
Since there exists a unique $\rho $ in $(0,1)$ which satisfies
the following equation:
$$\lambda_{0}(M)- \frac{48M}{\pi^{2}}s_{0}\arctan\rho-\frac{2T(M)\rho}{(1-\rho)^{3}}=0,
$$
we see that $H(z_{1})\neq H(z_{2})$ and so, $H(z)$ is
univalent for $|z|<\rho_{0}.$

Furthermore, we observe  that for any $z$ with $|z|=\rho_{0},$
\begin{eqnarray*}
|H(z)|&=&\rho_{0}^{2(p-1)}\big |zG_{z}(0)-\overline{z}G_{\overline{z}}(0)+
z(G_{z}(z)-G_{z}(0))-\overline{z}(G_{\overline{z}}(z)-G_{\overline{z}}(0))\big |\\
&\geq&\rho_{0}^{2p-1}\Big
[\lambda_{0}(M)-\frac{16M}{\pi^{2}}s_{0}\arctan\rho_{0} \Big ]\\
&=& R.
\end{eqnarray*}
The proof of the theorem is complete.
\epf
\begin{table}
\center{
\begin{tabular}{|c|c||c|c|c|c|}
\hline
$M$ & $p$ &$\rho=\rho(M,p)$ & $R=R(M,\rho (M,p))$  & $\rho '$ & $R'$
\\
\hline
1.1296& 2& 0.0281673 & 0.0000106985 & 0.0194864  & 3.54498$\times 10^{-6}$
\\
2& 2& 0.00856025  & 1.73218$\times 10^{-7}$ & 0.00623202 & 6.5415$\times 10^{-8}$
\\
2.2976& 2& 0.00646284 & 6.4986$\times 10^{-8}$ & 0.0047235 & 2.47902$\times 10^{-8}$
\\
3& 2&0.0037942 & 1.00669$\times 10^{-8}$  & 0.00277162 & 3.83502$\times 10^{-9}$
\\
1.1296& 3& 0.0281673    & 8.48819$\times 10^{-9}$ & -- & --
\\
2& 3& 0.00856025  & 1.2693$\times 10^{-11}$ & -- & --
\\
2.2976& 3& 0.00646284  & 2.71435$\times 10^{-12}$ & -- & --
\\
3&3 &0.0037942 & 1.44922$\times 10^{-13}$ & -- & --
\\
\hline
\end{tabular}
}
\bigskip
\caption{Values of $\rho$ and $R$ for Theorem \ref{thm4} for $p=2,3$,
and the corresponding values of $\rho '$ and $R '$ of \cite[Theorem 1.2]{CPW0} (for $p=2$)
\label{table2}
}
\end{table}

We remark that Theorem \ref{thm4} is an improved version of
\cite[Theorem 1.2]{CPW0} when $p=2$. In order to be more explicit, we refer to Table \ref{table2}
in which the third and fourth columns refer to values obtained from Theorem \ref{thm4} for cases $p=2,3$
for certain choices of $M$, while the right two columns correspond to the values obtained
from \cite[Theorem 1.2]{CPW0} for the case $p=2$.

\section{The Region of Variability}\label{sec-4}

\begin{defn}
Let $\mathcal{H}_{p}$ denote the set of all $p$-harmonic mappings of
the unit disk $\mathbb{D}$ with the normalization
$f_{z^{p-1}}(0)=(p-1)!$ and $|f(z)|\leq 1$ for $|z|<1$. Here we
prescribe that $\mathcal{H}_{0}=\emptyset$.

For a fixed point $z_{0}\in\mathbb{D}$, let
$$V_{p}(z_{0})=\{f(z_{0}):\ f\in \mathcal{H}_{p}\setminus\mathcal{H}_{p-1}\}.
$$
\end{defn}

Now, we have

\begin{thm}\label{thm3.1}
\begin{enumerate}
\item[{\rm \textbf{(a)}}] If $p=1$, then $V_{1}(z_{0})=\{1\};$
\item[{\rm \textbf{(b)}}] If  $p\geq 2$, $V_{p}(z_{0})=\mathbb{\overline{D}}$.
\end{enumerate}
\end{thm}
\bpf We first prove \textbf{(a)}.  Let $f\in\mathcal{H}_{1}$ and
$f(z)=\sum_{n=0}^{\infty}a_{n}z^{n}+\sum_{n=1}^{\infty}\overline{b}_{n}\overline{z}^{n}$.
By Parseval's identity and the hypotheses $|f(z)|\leq 1$ and
$f(0)=1$, we have
\begin{eqnarray*}
\lim_{r\rightarrow
1-}\frac{1}{2\pi}\int_{0}^{2\pi}|f(re^{i\theta})|^{2}\,d\theta
&=&\lim_{r\rightarrow 1-}\frac{1}{2\pi}\int_{0}^{2\pi}\big
(|h(re^{i\theta})|^{2}+|g(re^{i\theta})|^{2} \big ) \,d\theta\\
&=&|a_{0}|^{2}+\sum_{n=1}^{\infty}\big (|a_{n}|^{2}+|b_{n}|^{2}\big
)\leq 1.
\end{eqnarray*}
This inequality implies that for any $n\geq 1$, $a_{n}=b_{n}= 0$
which gives that $f(z)\equiv 1$ for  $z\in\mathbb{D}$. Thus, we have
$V_{1}(z_{0})=\{1\}$.

In order to prove \textbf{(b)}, we consider the function
\begin{eqnarray*}
\phi (z)=\frac{z^{p-1}-w}{1-w\overline{z}^{p-1}}
=|z|^{2(p-1)}\sum_{n=1}^{\infty}w^{n}\overline{z}^{(n-1)(p-1)}+z^{p-1}-w-\sum_{n=1}^{\infty}w^{n+1}\overline{z}^{(p-1)n},
\end{eqnarray*}
where $w\in\mathbb{\overline{D}}$ and $p\geq2$.

Then $\phi_{z^{p-1}}(0)=(p-1)!$,  $\Delta ^p \phi=0$ and therefore,
$\phi \in \mathcal{H}_{p}\setminus\mathcal{H}_{p-1}.$ For each fixed
$a\in\mathbb{\overline{D}}$, $z\mapsto
f_a(z)=(z^{p-1}-a)/(1-a\overline{z}^{p-1})$ is a $p$-harmonic
mapping and $f_a(\mathbb{D})\subset\mathbb{D}$.

 Obviously, $a\mapsto f_a(z_{0})
=\frac{z_{0}^{p-1}-a}{1-a\overline{z_{0}}^{p-1}}$ is a conformal
automorphism of $\mathbb{D}$ and the image of
$\overline{\mathbb{D}}$ under $f_a(z_{0})$ is
$\overline{\mathbb{D}}$ itself. By hypotheses, we obtain that for
any $g\in\mathcal{H}_{p}\setminus\mathcal{H}_{p-1}$,
$g(z_{0})\in\mathbb{\overline{D}}$. Hence $V_0(z_{0})$ coincides
with $\mathbb{\overline{D}}$. The proof of this theorem is complete.
\epf

By the method of proof used in Theorem \ref{thm3.1}\textbf{(a)}, we
obtain the following generalization of Cartan's uniqueness theorem
(see
 \cite{Ch} or \cite[p.~23]{R}) for harmonic mappings.

\begin{thm}\label{thm3.1x}
Let $f$ be a harmonic mapping in $\mathbb{D}$ with
$f(\mathbb{D})\subseteq \mathbb{D}$ and $f_{z}(0)=1$. Then $f(z)=z$
in $\mathbb{D}$.
\end{thm}

\section{Estimates for Bloch norm for bi- and tri-harmonic mappings}\label{sec6}

In the case of $p$-harmonic Bloch mappings, the authors in
\cite{CPW2} obtained the following result.

\begin{thm}\label{thm1.0}
Let $f$ be a $p$-harmonic mapping in $\mathbb{D}$ of the form {\rm
(\ref{eq2.0x})} satisfying $B_{f}<\infty$, where
$$B_{f}:=\sup_{z,w\in\mathbb{D},\ z\neq w}\frac{|f(z)-f(w)|}{\rho(z,w)}<\infty ~\mbox{ with }~
\rho(z,w)=\frac{1}{2}\log\left(\frac{1+|\frac{z-w}{1-\overline{z}w}|}
{1-|\frac{z-w}{1-\overline{z}w}|}\right).
$$
Then \beq\label{eq2}
B_f&:=&\sup_{z\in\mathbb{D}}(1-|z|^{2})\left \{\left |\sum_{k=1}^{p}|z|^{2(k-1)}(G_{p-k+1})_{z}(z) \right .\right .  \nonumber\\
&& \left . +
\sum_{k=1}^{p}(k-1)\overline{z}|z|^{2(k-2)}G_{p-k+1}(z)\right |
+\left |\sum_{k=1}^{p}|z|^{2(k-1)}(G_{p-k+1})_{\overline{z}}(z)
\right . \nonumber \\ \nonumber
&& + \left . \left .\sum_{k=1}^{p}(k-1)z|z|^{2(k-2)}G_{p-k+1}(z)\right |\right \} \\
&\geq& \sup_{z\in\mathbb{D}}(1-|z|^{2})\left | \Big
|\sum_{k=1}^{p}|z|^{2(k-1)}(G_{p-k+1})_{z}(z)\Big | -
 \Big |\sum_{k=1}^{p}|z|^{2(k-1)}(G_{p-k+1})_{\overline{z}}
 (z)\Big |\right|
\eeq and {\rm (\ref{eq2})} is sharp. The  equality sign in {\rm
(\ref{eq2})} occurs when $f$ is analytic or anti-analytic.

Furthermore, if for each $k\in \{1, 2, \ldots , p\}$,  the harmonic
functions $G_{p-k+1}$ in {\rm (\ref{eq2.0x})} are such that
$|G_{p-k+1}(z)|\leq M$, then \be\label{eq3} B_{f}\leq 2M
\phi_p(y_0). \ee Here $y_0$ is the unique root in $(0,1)$ of the
equation $\phi'_p(y)=0$, where \be\label{eq-phi1} \phi_p(y) =
\frac{2}{\pi}\sum_{k=1}^p y^{2(k-1)}+y(1-y^2)\sum_{k=2}^p
(k-1)y^{2(k-2)}. \ee
The bound in {\rm (\ref{eq3})} is sharp when $p=1$, where $M$ is a
positive constant. The extremal functions are
$$f(z)=\frac{2M\alpha}{\pi}{\rm Im}\left(\log\frac{1+S(z)}{1-S(z)}\right),
$$
where $|\alpha|=1$ and $S(z)$ is a conformal automorphism of
$\mathbb{D}$.
\end{thm}

In order to emphasize the importance of this result, we recall that,
when $p=1$, (\ref{eq2}) (resp. (\ref{eq3})) is a generalization of
\cite[Theorem 1]{Co} (resp. \cite[Theorem 3]{Co}). In the case of
$p=2$ of Theorem \ref{thm1.0}, after some computation, one has the
following simple formulation for biharmonic mappings.

\bcor\label{cor1a}
Let $f=H+|z|^2G$ be a biharmonic mapping of $\ID$
such that $B_{f}<\infty$. Then, we have \be\label{extra-eq2} B_f
\geq  \sup_{z\in\mathbb{D}}(1-|z|^{2})\left | |H_z + |z|^2G_z| -
|H_{\overline{z}} + |z|^2G_{\overline{z}}| \right | \ee and
\be\label{extra-eq3} B_{f}\leq
\frac{4M}{27\pi^3}\left(8+36\pi^2+\left(4+3\pi^2\right)^{3/2}\right)
\approx 30.7682 M. \ee \ecor\bpf According to our notation,
(\ref{eq2}) is equivalent to (\ref{extra-eq2}). In order to prove
(\ref{extra-eq3}), we first observe that (\ref{eq3}) is equivalent
to
$$B_{f}\leq 2M \sup_{0<y<1} \phi_2(y),
$$
where
$$\phi_2(y)=\frac{2}{\pi}(1+y^2)+y(1-y^2).
$$
Now, to find $\ds\sup_{0<y<1} \phi_2(y)$, we compute the derivative
$$\phi_2'(y) = 1+\frac{4}{\pi}y-3y^2 = -3\left(y-y_0\right)\left(y-\frac{2-\sqrt{4+3\pi^2}}{3\pi}\right)
$$
so that $\phi_2'(y)\geq 0$ for $0\leq y\leq y_0$ and $\phi_2'(y)\leq
0$ for $y_0\leq y<1$. Hence
$$y_0=\frac{2+\sqrt{4+3\pi^2}}{3\pi}\approx 0.82732
$$
is the critical point of $\phi_2(y)$. Consequently, $\phi_2(y)\leq
\phi_2(y_0).$ A simple calculation shows that \beqq
\phi_2(y_0)&=&\frac{2}{\pi}(1+y_0^2)+y_0(1-y_0^2)\\
&=&
\frac{2}{\pi}\left(\frac{8+12\pi^2+4\sqrt{4+3\pi^2}}{9\pi^2}\right)+
\left(\frac{2}{3\pi}+\frac{\sqrt{4+3\pi^2}}{3\pi}\right)\left(\frac{6\pi^2-8-4\sqrt{4+3\pi^2}}{9\pi^2}\right)\\
&=& \frac{2}{27\pi^3}\left(16+42\pi^2+8\sqrt{4+3\pi^2}+\sqrt{4+3\pi^2}\left(3\pi^2-4-2\sqrt{3\pi^2+4}\right)\right)\\
&=&\frac{2}{27\pi^3}\left(8+36\pi^2+\left(4+3\pi^2\right)^{3/2}\right)
\approx 15.3841 \eeqq and therefore, $B_{f}\leq 2M\phi_2(y_0)$ which
is the desired inequality (\ref{extra-eq3}). The result follows.
\epf

In the case of $p=3$ of Theorem \ref{thm1.0}, we have

\bcor\label{cor2a}
 Let $f=H+|z|^2G+|z|^4K$ be a triharmonic $($i.e.
$3$-harmonic$)$ mapping of the unit disk $\ID$ such that
$B_{f}<\infty$, where $H$, $G$ and $K$ are harmonic in $\ID$. Then
we have \be\label{extra-eq4} B_f \geq
\sup_{z\in\mathbb{D}}(1-|z|^{2})\left | |H_z + |z|^2G_z
+|z|^4K_z\right | - \left |H_{\overline{z}}+|z|^2G_{\overline{z}}
+|z|^4K_{\overline{z}}| \right | \ee and \be\label{extra-eq5}
B_{f}\leq 2M \phi_3(y_1) \approx 4.037006 M, \ee where $\phi_3(y_1)
= \sup_{0<y<1} \phi_3(y)$ and
$$\phi_3(y)=\frac{2}{\pi}(1+y^2+y^4)+y(1+y^2-2y^4).
$$
\ecor\bpf Set $p=3$ in Theorem \ref{thm1.0}. Then, (\ref{extra-eq4})
is equivalent to (\ref{eq2}) and therefore, it suffices to prove
(\ref{extra-eq5}). The choice $p=3$ in (\ref{eq3}) shows that
$$B_{f}\leq 2M \sup_{0<y<1} \phi_3(y),
$$
where $\phi_3(y)$ is obtained from {\rm (\ref{eq-phi1})}.

We see that $\phi_3(y)$ has a unique positive root in $(0,1)$. Also,
$$\phi_3'(y) = \frac{4}{\pi}(y+2y^3)+1+3y^2-10y^4.
$$
Computations show that    $\phi_3'(y)\geq 0$ for $0\leq y\leq y_1$
and $\phi_3'(y)\leq 0$ for $y_1\leq y<1$. Hence

$$y_1\approx 0.891951$$
is the only
critical point of $\phi_3(y)$ in the interval $(0,1)$. It follows that
$$\phi_3(y)\leq \phi_3(y_1)\approx 2.018503.
$$
Thus, $B_{f}\leq 2M\phi_3(y_1)$ which is the desired inequality
(\ref{extra-eq5}). \epf

\end{document}